\documentclass[letterpaper, 10 pt, conference]{ieeeconf}  

\IEEEoverridecommandlockouts                              
\usepackage{graphics} 
\usepackage{epsfig} 
\usepackage{mathptmx} 
\usepackage{times} 
\usepackage{amsmath} 
\usepackage{amssymb}  
\usepackage{placeins} 

\usepackage{xcolor}
\usepackage{siunitx} 
\usepackage[ruled,linesnumbered]{algorithm2e} 

\usepackage{svg}
\graphicspath{{../img/}}

\let\mathcal\relax\usepackage[cal=cm]{mathalfa}

\usepackage{bm}

\usepackage{url}

\usepackage[multiple]{footmisc}

\usepackage{doi}

\newif\ifcommentandrea
\commentandreatrue

\newtheorem{remark}{Remark}
\newtheorem{lemma}{Lemma}

\newtheorem{assumption}{Assumption}

\title{\LARGE \bf
Zero-Order Optimization for \\Gaussian Process-based Model Predictive Control
}

\author{A. Lahr$^1$, A. Zanelli$^1$, A. Carron$^1$ and M. N. Zeilinger$^1$ 
\thanks{$^1$Institute for Dynamic Systems and Control, ETH Zürich, Zürich \mbox{CH-8092}, Switzerland. E-mail: {\tt\small amlahr@ethz.ch}.}%
\thanks{This project has received funding by the European Union's Horizon 2020 research and innovation programme, Marie~Sk\l{}odowska-Curie grant agreement No. 953348,~\mbox{ELO-X}.}%
}

\begin{document}

\maketitle
\thispagestyle{empty}
\pagestyle{empty}

\begin{abstract}

    By enabling constraint-aware online model adaptation, model predictive control using Gaussian process~(GP) regression has exhibited impressive performance in real-world applications and received considerable attention in the learning-based control community. Yet, solving the resulting optimal control problem in real-time generally remains a major challenge, 
    due to i)~the increased number of augmented states in the optimization problem, as well as ii)~computationally expensive evaluations of the posterior mean and covariance and their respective derivatives. To tackle these challenges, we employ i)~a~tailored Jacobian approximation in a sequential quadratic programming (SQP) approach and combine it with ii)~a~parallelizable GP inference and automatic differentiation framework. Reducing the numerical complexity with respect to the state dimension $n_x$ for each SQP iteration from $\mathcal{O}(n_x^6)$ to $\mathcal{O}(n_x^3)$,  and accelerating GP evaluations on a graphical processing unit, the proposed algorithm computes suboptimal, yet feasible, solutions at drastically reduced computation times and exhibits favorable local convergence properties. Numerical experiments verify the scaling properties and investigate the runtime distribution across different parts of the algorithm.  

\end{abstract}

\section{Introduction}

Real-world applications of model predictive control using Gaussian processes~(GP-MPC), such as vision-based robot path-tracking~\cite{ostafewLearningbasedNonlinearModel2014}, trajectory tracking using a robotic arm~\cite{carronDataDrivenModelPredictive2019}, autonomous racing~\cite{kabzanLearningBasedModelPredictive2019,kabzanAMZDriverlessFull2020}, or high-speed quadrotor flight~\cite{torrenteDataDrivenMPCQuadrotors2021}, have showcased its potential to leverage closed-loop data for constraint-aware online model adaptation. 

Yet, the computational cost associated with GP inference remains a major challenge. To fully exploit the rich, state-dependent uncertainty description induced by the posterior covariance of the GP, it should be included in the optimal control problem (OCP) formulation and propagated through the dynamics model over the prediction horizon~\cite{hewingCautiousModelPredictive2020}. 
However, incorporating the covariance propagation into the OCP constitutes a major limitation for real-time implementation of the algorithm on embedded hardware, both in terms of additional optimization variables capturing the state covariance, as well as computationally expensive evaluation and differentiation of the GP~posterior covariance. 

As a result, practical implementations resort to various heuristics to speed up the algorithm. Popular approaches include GP~approximations with a fixed number of basis functions%
~\cite{frohlichModelLearningContextual2022, vaskovFrictionAdaptiveStochasticPredictive2022}, 
reduction of optimization variables by fixing the state covariances in the optimization problem based on their predicted value at the last MPC iteration~\cite{kabzanLearningBasedModelPredictive2019, hewingCautiousModelPredictive2020}, or to completely ignore the uncertainty description provided by the GP posterior covariance in the control algorithm~\cite{ostafewLearningbasedNonlinearModel2014},~\cite{torrenteDataDrivenMPCQuadrotors2021}.

\subsection{Contributions}

This paper addresses the aforementioned challenges by 
\begin{itemize}
    \item[i)] applying a tailored inexact sequential quadratic programming (SQP) algorithm~\cite{fengInexactAdjointbasedSQP2020,zanelliZeroOrderRobustNonlinear2021} 
    to solve the GP-MPC optimal control problems, returning suboptimal, yet feasible solutions at convergence,
    \item[ii)] speeding up GP computations by efficient inference and automatic differentiation~(AD)~\cite{paszkePyTorchImperativeStyle2019,gardnerGPyTorchBlackboxMatrixMatrix2018}, as well as parallelization on a graphical processing unit (GPU),
    \item[iii)] showing that the resulting algorithm maintains favorable local convergence properties if the GP posterior covariance and process noise covariance at the solution are sufficiently small.
\end{itemize}


\subsection{Related Work}


A tailored Jacobian approximation for stochastic and robust nonlinear model predictive control (NMPC) has been originally developed in~\cite{fengInexactAdjointbasedSQP2020,zanelliZeroOrderRobustNonlinear2021} where it is referred to as a zero-order method. Based on this approach, \cite{messererEfficientAlgorithmTubebased2021} presents a robust NMPC algorithm that improves the ellipsoidal uncertainty propagation by efficiently optimizing for optimal linear feedback policies used in the predictions. In the stochastic setting, in~\cite{quirynenUncertaintyPropagationLinear2021} the tailored Jacobian approximation is coupled with a more accurate uncertainty propagation based on linear-regression Kalman filtering, showing real-time feasible timings for an automotive example. With respect to the aforementioned works, this paper presents a direct extension of~\cite{fengInexactAdjointbasedSQP2020, zanelliZeroOrderRobustNonlinear2021} to the setting of Gaussian process-based MPC. As such, the improvements presented in~\cite{messererEfficientAlgorithmTubebased2021, quirynenUncertaintyPropagationLinear2021} can be straightforwardly applied to the results in this paper as well.

To the best of the authors' knowledge, \cite{vaskovFrictionAdaptiveStochasticPredictive2022}~is the only publication that applies the tailored Jacobian approximation to a variant of GP-MPC. Therein, the tailored stochastic NMPC algorithm is used in conjunction with an approximate GP based on a finite set of basis functions for an autonomous driving application, where the distributions of the weights are estimated jointly with the system state using a particle filter. Conceptually, the presented work differs from~\cite{vaskovFrictionAdaptiveStochasticPredictive2022} insofar as we apply exact GP inference with an approximate uncertainty propagation, while in~\cite{vaskovFrictionAdaptiveStochasticPredictive2022}, a finite-dimensional GP approximation is used and the uncertainty propagation is performed using a particle filter. Moreover, we additionally focus on a theoretical analysis of the convergence properties and provide a high-performance implementation of the proposed algorithm. 
\section{Problem formulation} \label{sec:problem}



We consider discrete-time, nonlinear dynamics of the form
\begin{align}
   x(k+1) = \psi(x(k), u(k)) + B\left(\eta(x(k), u(k)) + w(k) \right),
\end{align}
where \mbox{$\psi: \mathbb{R}^{n_x} \times \mathbb{R}^{n_u} \rightarrow \mathbb{R}^{n_x}$}, \mbox{$\eta: \mathbb{R}^{n_x} \times \mathbb{R}^{n_u} \rightarrow \mathbb{R}^{n_w}$} describe the known and unknown parts of the system dynamics as functions of the state~\mbox{$x(k) \in \mathbb{R}^{n_x}$} and input~\mbox{$u(k) \in \mathbb{R}^{n_u}$} at time step $k$, respectively, and are assumed to be twice continuously differentiable. The process noise $w(k) \sim \mathcal{N}(0,\Sigma^w)$ is assumed to be i.i.d. in each component, i.e., with diagonal covariance matrix~\mbox{$\Sigma^w = \mathrm{diag}\left( \sigma_1^2, \ldots, \sigma_{n_w}^2 \right)$}; it affects the states via the matrix~\mbox{$B \in \mathbb{R}^{n_x \times n_w}$}, allowing to model the process noise in a lower-dimensional space. The system is subject to $n_h$ individual chance constraints, i.e.,
\begin{align}
    \mathrm{Pr}(h_j(x(k),u(k)) \leq 0) \geq p_j, \label{eq:ChanceConstraints}
\end{align}
for all constraints $j=1,\ldots,n_h$ and time steps $k$, with satisfaction probability~\mbox{$0 < p_j \leq 1$}. This formalism also captures hard constraints, for example on the inputs, by setting~\mbox{$p_j = 1$} for the respective constraint.

The key idea of GP-MPC is to model the unknown dynamics $\eta$ as a Gaussian process
\begin{align}
    d(x,u) \sim \mathcal{GP} \left(\mu^d(x,u),\Sigma^d(x,u) \right),
\end{align}
with posterior mean $\mu^d: \left(\mathbb{R}^{n_x} \times \mathbb{R}^{n_u} \right) \rightarrow \mathbb{R}^{n_w}$ and posterior covariance~\mbox{$\Sigma^d: (\mathbb{R}^{n_x} \times \mathbb{R}^{n_u}) \rightarrow \mathbb{R}^{n_w \times n_w}$}, respectively. Note that by~$d(x,u)$ we denote a GP conditioned on data points already; see e.g.~\cite[Chap.~2]{rasmussenGaussianProcessesMachine2006} for the inference formulas.
At every sampling time, the problem to be solved can be formulated as a stochastic OCP\footnote{In~\eqref{eq:SOCP}, terminal constraints have been omitted for simplicity.}, 
where a sequence of control polices~\mbox{$\bm{\pi} := \{ \pi_i \}_{i=0}^{N-1}$}, \mbox{$\pi: \mathbb{R}^{n_x} \rightarrow \mathbb{R}^{n_u}$}, 
is determined 
over a given set of admissible policies, that minimizes the expected value of a user-defined cost over a finite time horizon of length $N$,%
\begin{subequations} \label{eq:SOCP}
    \begin{align}
        \underset{\bm{\pi}}{\min} \quad & \mathbb{E} \Bigg[ c_f( x_N) + \sum_{i=0}^{N-1} c_i(x_i, u_i )  \Bigg] \\
            \mathrm{s.t.} \quad &\forall i \in \{0,\ldots,N-1\}, \\
            & u_i = \pi_i(x_i), \\
            & x_{i+1} = \psi(x_i, u_i) + B \left( d(x_i,u_i) + w_i \right), \label{eq:SOCP_dynamics} \\
            & \mathrm{Pr}(h_j(x_i,u_i) \leq 0) \geq p_j, \> j=1,\ldots,n_h, \label{eq:SOCP_Chanceconstraints} \\
            & x_0 = x(k).
    \end{align}
\end{subequations}


As computing exact solutions to the stochastic OCP~\eqref{eq:SOCP} is generally intractable, in the following we present a common approach to approximate~\eqref{eq:SOCP} by a deterministic problem, as previously proposed by~\cite{hewingCautiousModelPredictive2020}. 

\subsection{Mean and covariance propagation}

Propagating the uncertainty introduced by the GP model through the nonlinear dynamics model in a computationally efficient way generally proves to be very challenging: As Gaussianity of the state distribution is lost after a nonlinear transformation, so is the possibility to capture the entire distribution by its first two moments. To retain computational tractability, we use a common linearization-based approximation of the propagation of the expected state and covariance~\cite{girardGaussianProcessPriors2002}, 
which results in the following deterministic update equations~\cite{hewingCautiousNMPCGaussian2018}, starting from \mbox{$x_0 \sim \mathcal{N}(\mu^x_0, \Sigma^x_0)$},%
\begin{align}
    \mu_{i+1}^x &= \psi(\mu_i^x, u_i) + B \mu^d(\mu_i^x, u_i), \\
	\Sigma_{i+1}^x & = \tilde{A}_i \Sigma_i^x \tilde{A}_i^\top + B \left( \Sigma^d(\mu_i^x, u_i) + \Sigma^w \right) B^\top , \label{eq:CovariancePropagation}
    \intertext{where}
    \tilde{A}_i &:= \left. \frac{\partial}{\partial x} \left( \psi(x,u_i) + B \mu^d(x, u_i) \right) \right|_{x=\mu^x_i}.
\end{align} 

\subsection{Chance constraints}

Given the covariance propagation above, it is possible to efficiently formulate the individual chance constraints~\eqref{eq:SOCP_Chanceconstraints} as deterministic constraints on the mean prediction by a suitable tightening, i.e., \mbox{$\bar{h}_j(\mu^x_i,u_i,\Sigma^x_i) \leq 0$}, where%
\begin{align*}
    \bar{h}_j(\mu^x_i,u_i,\Sigma^x_i) := h_j(\mu^x_i, u_i) + \alpha_j \sqrt{C_j(\mu^x_i, u_i) \Sigma^x_i C_j(\mu^x_i, u_i)^\top}
\end{align*}
and $C_j(\mu^x_i, u_i) := \frac{\partial h_j}{\partial x} (\mu^x_i, u_i)$. For general probability distributions, the tightening factor 
\mbox{$\alpha := \sqrt{\frac{p_j}{1 - p_j}}$} 
can be chosen based on the Chebyshev inequality; for Gaussian distributions, setting \mbox{$\alpha_j :=  \Phi^{-1} (p_j)$}, where $\Phi^{-1}(\cdot)$ is the inverse cumulative density function of a standard normal Gaussian variable, is a less conservative choice. 



\subsection{Optimal control problem}

With the above simplifications and an approximate formulation of the expected cost, we arrive at the following deterministic approximation of the stochastic OCP~\eqref{eq:SOCP} in terms of the predicted state mean \mbox{$\bm{\mu} := \{ \mu^x_i \}_{i=0}^N$}, state covariance \mbox{$\bm{\Sigma} := \{ \Sigma^x_i \}_{i=0}^N$} and control input sequence~\mbox{$\bm{u} := \{u_i\}_{i=0}^{N-1}$}, cf.~\cite{hewingCautiousModelPredictive2020},%
\begin{subequations} \label{eq:OCP_GPMPC}
    \begin{align}
        \underset{\bm{\mu},\bm{\Sigma},\bm{u}}{\min} \quad & c_f(\mu_N^x) + \sum_{i=0}^{N-1} c_i(\mu_i^x, u_i)  \\
        \mathrm{s.t.} \quad & \forall i \in \{0,\ldots,N-1\}, \\
        & \mu_{i+1}^x = \psi(\mu_i^x, u_i) + B \mu^d(\mu_i^x, u_i), \label{eq:OCP_GPMPC_mean} \\
        & \Sigma_{i+1}^x = \tilde{A}_i \Sigma_i^x \tilde{A}_i^\top + B \left( \Sigma^d(\mu_i^x, u_i) + \Sigma^w \right) B^\top \label{eq:OCP_GPMPC_Sigma}, \\
        & \bar{h}_j(\mu^x_i,u_i,\Sigma^x_i) \leq 0, \> j=1,\ldots,n_h, \\
        & \mu_0^x = x(k) \label{eq:OCP_GPMPC_mean_0}, \\
        & \Sigma_0^x = 0 \label{eq:OCP_GPMPC_Sigma_0}.
    \end{align}
\end{subequations}

\begin{remark}
    While optimizing over control input sequences in~\eqref{eq:OCP_GPMPC} has computational advantages compared to optimization over a sequence of feedback policies in~\eqref{eq:SOCP}, it can lead to underconfident predictions and an overly conservative controller. 
    As a remedy, linear state feedback can be incorporated into the controller~\cite{hewingCautiousModelPredictive2020}. In that case, hard input constraints can generally not be satisfied and should be replaced with individual chance constraints; the adaptation of the OCPs~\eqref{eq:SOCP} and~\eqref{eq:OCP_GPMPC} is straightforward and can be found in~\cite{hewingCautiousModelPredictive2020}. 
\end{remark}

Solving~\eqref{eq:OCP_GPMPC} in a receding horizon fashion yields a highly performant and adaptive, yet uncertainty-aware, control strategy, demonstrated by real-world applications such as~\cite{kabzanLearningBasedModelPredictive2019, torrenteDataDrivenMPCQuadrotors2021}. 
Nevertheless, for a real-time GP-MPC implementation, there remain major computational challenges to be addressed.

First, the OCP~\eqref{eq:OCP_GPMPC} 
has $n_x + n_u + (n_x + n_x^2)/2$ optimization variables at each prediction stage, structurally equivalent to OCPs arising from a direct multiple shooting discretization~\cite{bockMultipleShootingAlgorithm1984}. For this set of problems, among the most competitive solvers at the time of this writing are interior point algorithms exploiting the block-sparse structure of the multiple shooting approach~\cite{kouzoupisRecentAdvancesQuadratic2018}.
Still, these algorithms' cubic computational complexity in the number of shooting nodes leads to a total computational complexity of $\mathcal{O}(n_x^6)$, which becomes prohibitive even for a moderate number of states. 

Second, even without updating the GP data points online, 
evaluating the GP posterior mean, covariance, and their respective Jacobians generally scales quadratically with the number of data points~\cite[p.~19]{rasmussenGaussianProcessesMachine2006}, posing strict limitations on the number of data points amenable for online inference. In this regard, solving~\eqref{eq:OCP_GPMPC} is particularly demanding, as computing the constraint Jacobians~\eqref{eq:OCP_GPMPC_Sigma} requires not only to compute the Jacobian of the nominal and GP dynamics but also, their respective second-order derivatives.

\begin{figure}
    \vspace{1mm}
    \centering
    \def\svgwidth{1.06\columnwidth}
\begingroup%
  \makeatletter%
  \providecommand\color[2][]{%
    \errmessage{(Inkscape) Color is used for the text in Inkscape, but the package 'color.sty' is not loaded}%
    \renewcommand\color[2][]{}%
  }%
  \providecommand\transparent[1]{%
    \errmessage{(Inkscape) Transparency is used (non-zero) for the text in Inkscape, but the package 'transparent.sty' is not loaded}%
    \renewcommand\transparent[1]{}%
  }%
  \providecommand\rotatebox[2]{#2}%
  \newcommand*\fsize{\dimexpr\f@size pt\relax}%
  \newcommand*\lineheight[1]{\fontsize{\fsize}{#1\fsize}\selectfont}%
  \ifx\svgwidth\undefined%
    \setlength{\unitlength}{281.91477168bp}%
    \ifx\svgscale\undefined%
      \relax%
    \else%
      \setlength{\unitlength}{\unitlength * \real{\svgscale}}%
    \fi%
  \else%
    \setlength{\unitlength}{\svgwidth}%
  \fi%
  \global\let\svgwidth\undefined%
  \global\let\svgscale\undefined%
  \makeatother%
  \begin{picture}(1,0.47936631)%
    \lineheight{1}%
    \setlength\tabcolsep{0pt}%
    \put(0,0){\includegraphics[width=\unitlength,page=1]{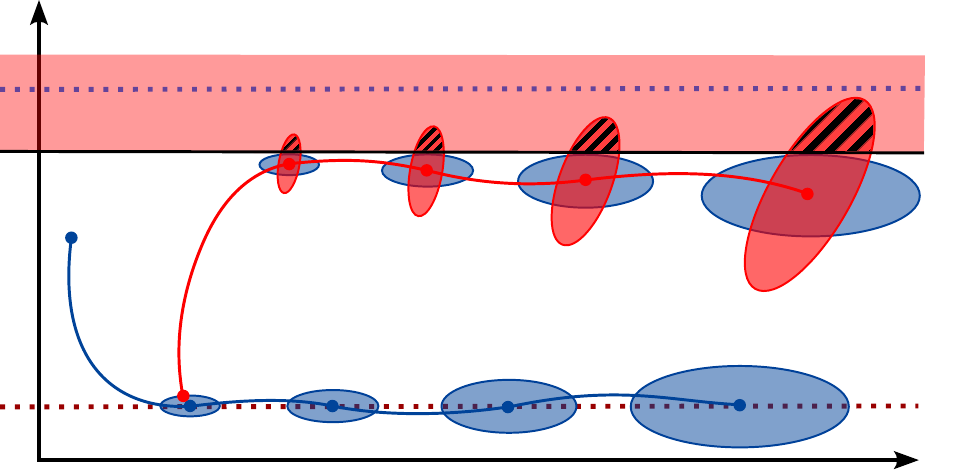}}%
    \put(0.05806389,0.45216157){\color[rgb]{0,0,0}\makebox(0,0)[lt]{\lineheight{1.25}\smash{\begin{tabular}[t]{l}$x_2$\end{tabular}}}}%
    \put(0.90686077,0.02940044){\color[rgb]{0,0,0}\makebox(0,0)[lt]{\lineheight{1.25}\smash{\begin{tabular}[t]{l}$x_1$\end{tabular}}}}%
  \end{picture}%
\endgroup%

    \caption{Infeasibility arising from fixing the covariances based on the previous MPC instance. The predicted state trajectory and covariances from the previous time step are drawn with solid lines in blue; predicted covariances around current predicted state trajectory in red; reference for previous/current time step with dashed lines; infeasible region in light red. When the linearized dynamics at the shooting nodes vary strongly from one time step to another, in this example based on a reference change, fixing the covariances based on the previous MPC instance might lead to significant prediction errors.}
    \label{fig:InfeasibilityHeuristic}
\end{figure}

A popular heuristic to ensure computational tractability despite the limitations stated above has been to propagate the state covariance~\eqref{eq:OCP_GPMPC_Sigma} outside the optimizer, based on input and state predictions from the last MPC iteration~\cite{kabzanLearningBasedModelPredictive2019},~\cite{hewingCautiousModelPredictive2020}. 
Fixing the covariance matrices $\Sigma^x_i$ at each stage in~\eqref{eq:OCP_GPMPC} eliminates the corresponding $(n_x + n_x^2)/2$ augmented states as well as the need to compute the equality constraint Jacobian derived from~\eqref{eq:OCP_GPMPC_Sigma}. However, this heuristic generally does not lead to feasible solutions for problem~\eqref{eq:OCP_GPMPC}, as Fig.~\ref{fig:InfeasibilityHeuristic} illustrates. 

To address the aforementioned challenges, 
we 
propose a tailored sequential quadratic programming~(SQP) algorithm that computes suboptimal, yet feasible, solutions to~\eqref{eq:OCP_GPMPC} with drastically reduced computational footprint.

\section{Zero-order method for GP-MPC} \label{sec:method}

In the following, we present a tailored optimization method for Gaussian process-based MPC based on inexact sequential quadratic programming. To avoid excessive notation in the following exposition, let us rewrite the GP-MPC optimal control problem~\eqref{eq:OCP_GPMPC} in the following compact form,%
\begin{subequations} \label{eq:OCP_compact}
    \begin{align}
        \underset{y, P}{\min} \quad  & c(y) \\ 
        \mathrm{s.t.} \quad & f(y) = 0, \label{eq:OCP_compact_equality_mean} \\
        & g(y,P) = 0, \label{eq:OCP_compact_equality_covariance} \\
        & \bar{h}(y,P) \leq 0. \label{eq:OCP_compact_inequality}
    \end{align}
\end{subequations}
Thereby we have concatenated all mean states and inputs into \mbox{$y := (y_0, \ldots, y_{N-1}, \mu^x_N)$} and vectorized covariances into \mbox{$P := (\overrightarrow{\Sigma^x_0}, \ldots, \overrightarrow{\Sigma^x_N})$}, where \mbox{$y_i := (\mu_i^x, u_i)$} for all \mbox{$i=0,\ldots,N-1$} and $\overrightarrow{(\cdot)}$ is the column-wise vectorization operator. 

\subsection{Sequential quadratic programming}

Sequential quadratic programming~\cite{boggsSequentialQuadraticProgramming1995} finds a solution to~\eqref{eq:OCP_compact} by iteratively solving quadratic programs that locally approximate the original NLP at the current solution estimate~$(\hat{y}, \hat{P})$. At each iteration, the solution~$(\Delta y, \Delta P)$ of the quadratic subproblem%
\begin{subequations}\label{eq:QP_standard}
    \begin{align}
       \underset{\Delta y, \Delta P}{\min} \quad & 
            \frac{1}{2}
            \Delta y^\top
            M_{yy}
            \Delta y
            +
            \left.
            \frac{\partial c}{\partial y}
            \right|_{y=\hat{y}}
            \Delta y\\ 
            \mathrm{s.t.} \> \quad &
            \begin{bmatrix}
                f  \\
                g 
            \end{bmatrix}
            (\hat{y},\hat{P})
            +
            \left.
            \begin{bmatrix}
                \frac{\partial f}{\partial y} & 0 \\
                \frac{\partial g}{\partial y} & \frac{\partial g}{\partial P}
            \end{bmatrix}
            \right|_{\overset{\scriptstyle y = \hat{y}}{P = \hat{P}}}
                \begin{bmatrix}
                    \Delta y \\
                    \Delta P
                \end{bmatrix}
            = 0, \label{eq:QP_equality_constraints} \\
            & 
            \bar{h}(\hat{y},\hat{P})
             + 
            \left.
            \begin{bmatrix}
                \frac{\partial \bar{h}}{\partial y} & \frac{\partial \bar{h}}{\partial P}
            \end{bmatrix}
            \right|_{\overset{\scriptstyle y = \hat{y}}{P = \hat{P}}}
            \begin{bmatrix}
                \Delta y \\
                \Delta P
            \end{bmatrix} 
            \leq 0,
    \end{align}
\end{subequations}
is then used to update the current solution estimate, where $M_{yy}$
denotes the chosen approximation of the Hessian of the Lagrangian of~\eqref{eq:OCP_compact}. 
The core SQP algorithm is summarized in Alg.~\ref{alg:sqp}.

\begin{algorithm} \label{alg:sqp}
    \caption{Standard SQP iteration}
    Initialize $\hat{y}, \hat{P}$\;
    \While{\text{termination criterion not met}}{
        Solve QP~\eqref{eq:QP_standard} for $\Delta y, \Delta P$\label{sqp:SolveQP}\;
        Set $y_+ := \hat{y} + \Delta y$, $P_+ := \hat{P} + \Delta P$\;
        Update $\hat{y} := y_+$, $\hat{P} := P_+$\;
    } 
\end{algorithm}


Under standard regularity assumptions, it can be shown that the SQP iteration locally converges to a Karush-Kuhn-Tucker (KKT) point of the original NLP~\eqref{eq:OCP_compact} if the Hessian approximation~$M_{yy}$ is sufficiently accurate, cf.~\cite[Thm.~3.2]{boggsSequentialQuadraticProgramming1995}.

Computationally, the SQP algorithm is ideally suited for NMPC, as the computational load can be reduced significantly by running only a fixed number of iterations before applying the control input to the plant, usually with negligible effect on the control performance, c.f. the Real-Time Iteration~\cite{bockConstrainedOptimalFeedback2007}. Still, as the number of optimization variables at each shooting node carries over from~\eqref{eq:OCP_GPMPC} to the quadratic subproblems in~\eqref{eq:QP_standard} when employing a sparsity-exploiting interior point solver, their cubic scaling with respect to the number of augmented states (or $\mathcal{O}(n_x^6)$ with respect to the state dimension) poses a major limitation on the nominal system dimension amenable for a GP-MPC application. 


\subsection{Zero-order algorithm}

To alleviate the limited scalability of optimizing over the state covariance matrices in GP-MPC, we employ an inexact SQP method that has been initially presented for stochastic~\cite{fengInexactAdjointbasedSQP2020} and robust NMPC~\cite{zanelliZeroOrderRobustNonlinear2021}. 
The core idea is that, by using a tailored Jacobian approximation in the SQP algorithm, the equality constraints in~\eqref{eq:QP_standard} corresponding to the covariance propagation~\eqref{eq:OCP_compact_equality_covariance} can be decoupled from the optimization problem. This allows for subsequent elimination of the associated optimization variables $\Delta P$ from the QP~\eqref{eq:QP_standard}, which we present in the following.

The tailored Jacobian approximation is obtained by neglecting the Jacobian of the covariance propagation~\eqref{eq:OCP_compact_equality_covariance} with respect to the nominal variables~$y$, i.e., by setting \mbox{$\frac{\partial g}{\partial y} := 0$} in%
~\eqref{eq:QP_equality_constraints}. This leads to a zero-order approximation of the corresponding equality constraints at each SQP iteration.
Under the approximate Jacobian, rearranging~\eqref{eq:QP_equality_constraints} yields an approximation of~$\Delta P$ based solely on the current linearization point $(\hat{y}, \hat{P})$, i.e.,
\begin{align}
    \Delta \tilde{P} = -\left( \left. \frac{\partial g}{\partial P} \right|_{\overset{\scriptstyle y=\hat{y}}{P=\hat{P}}} \right)^{-1} g(\hat{y},\hat{P}), \label{eq:DeltaP_InverseConstraints}
\end{align}
which can be used to eliminate $\Delta P \approx \Delta \tilde{P}$ from the QP subproblems~\eqref{eq:QP_standard}. 

Instead of solving the linear system in~\eqref{eq:DeltaP_InverseConstraints}, however, let us show how~$\Delta \tilde{P}$ may be obtained more efficiently and intuitively by making use of the fact that $g(y,P)$ corresponds to a vectorized version of the covariance propagation~\eqref{eq:CovariancePropagation}. 
We start by noting that the corresponding, vectorized equality constraints~\eqref{eq:OCP_compact_equality_covariance} may be written as%
\begin{align}
    g(y,P) = A(y) P + b(y) = 0, \label{eq:prop_linear}
\end{align}
where $A(y) \in \mathbb{R}^{n_x^2 (N+1) \times n_x^2 (N+1)}$ is invertible and contains the system dynamics and $b(y) \in \mathbb{R}^{n_x^2 (N+1)}$, the vectorized process noise and GP posterior covariances for all prediction steps; see Appendix~\ref{sec:appendix_matrices} %
for an explicit construction. Due to linearity of $g(y,P)$ in $P$, obtaining $\Delta \tilde{P}$ in~\eqref{eq:DeltaP_InverseConstraints} based on the current linearization point~$(\hat{y},\hat{P})$ is equivalent to solving for%
\begin{align}
    \tilde{P}_+ = -A(\hat{y})^{-1} b(\hat{y}) \label{eq:CovariancePropagation_Inverse}
\end{align} 
and setting $\Delta \tilde{P} := \tilde{P}_+ - \hat{P}$. Since~\eqref{eq:prop_linear} is obtained by vectorization and stacking of~\eqref{eq:CovariancePropagation} for all stages, it can easily be verified that solving for $\tilde{P}_+$ corresponds to propagating the predicted covariance as given by equation~\eqref{eq:CovariancePropagation}, based on the current iterate $\hat{y}$ and initial covariance \mbox{$\Sigma^x_0 := 0$}, followed by an application of the vectorization operator. 

After eliminating 
$\Delta P \approx \Delta \tilde{P}$, problem~\eqref{eq:QP_standard} reduces to%
\begin{subequations}\label{eq:QP_zoro}
    \begin{align}
        \underset{\Delta y}{\min} \quad &
        \frac{1}{2}
        \Delta y^\top
        M_{yy}
        \Delta_y
        +
        \left.
            \frac{\partial c}{\partial y}
            \right|_{y=\hat{y}}
        \Delta y \\
        \mathrm{s.t.} \quad &
        f(\hat{y}) + \left. \frac{\partial f}{\partial y} \right|_{\overset{\scriptstyle y=\hat{y}}{P=\hat{P}}} \Delta y 
        = 0, \\
        & \bar{h}(\hat{y}, \hat{P}) +
        \left.
        \begin{bmatrix}
            \frac{\partial \bar{h}}{\partial y} & \frac{\partial \bar{h}}{\partial P}
        \end{bmatrix}
        \right|_{\overset{\scriptstyle y=\hat{y}}{P=\hat{P}}}
        \begin{bmatrix}
            \Delta y \\
            \Delta \tilde{P}
        \end{bmatrix} \leq 0. \label{eq:QP_zoro_ineq}
    \end{align}
\end{subequations}
The modified iteration of the SQP algorithm alternates the solution of QP~\eqref{eq:QP_zoro} with the covariance propagation based on~\eqref{eq:CovariancePropagation} until convergence, as summarized in Alg.~\ref{alg:gpzoro}. Due to the reduced size of the QPs, we recover again the computational complexity of $\mathcal{O}(n_x^3)$ when employing sparsity-exploiting interior point QP solvers. This provides a drastic improvement compared with the complexity of $\mathcal{O}(n_x^6)$ for the QPs in the original problem~\eqref{eq:QP_standard}, as also demonstrated by the numerical examples in Section~\ref{sec:experiments}. 

\begin{algorithm} \label{alg:gpzoro}
    \caption{Modified SQP iteration}
    Initialize $\hat{y}, \hat{P}$\;
    \While{\text{termination criterion not met}}{
        Obtain $\tilde{P}_+ = A(\hat{y})^{-1} b(\hat{y})$ from~\eqref{eq:CovariancePropagation}\label{sqp:PropagateP}\;
        Set $\Delta \tilde{P} := \tilde{P}_+ - \hat{P}$\;
        Solve QP~\eqref{eq:QP_zoro} for $\Delta y$\label{sqp:SolveQP_zoro}\;
        Set $y_+ := \hat{y} + \Delta y$\;
        Update $\hat{y} := y_+$, $\hat{P} := \tilde{P}_+$\;
    }
\end{algorithm}



\section{Local convergence properties} \label{sec:convergence}

At convergence, the solution obtained from Alg.~\ref{alg:gpzoro} is guaranteed to be feasible and suboptimal for the original NLP~\eqref{eq:OCP_GPMPC}~\cite{zanelliZeroOrderRobustNonlinear2021},~\cite{bockConstrainedOptimalFeedback2007}. This is in contrast to existing heuristics discussed in Fig.~\ref{fig:InfeasibilityHeuristic}, where feasibility of the obtained solution with respect to~\eqref{eq:OCP_GPMPC} cannot be guaranteed.  

Whether Alg.~\ref{alg:gpzoro} converges to a feasible point of~\eqref{eq:OCP_GPMPC}, however, depends on the error in the Jacobian approximation. 
In~\cite{zanelliZeroOrderRobustNonlinear2021}, this has been studied for equality constraints of the form \mbox{$0 = A(y) P - \tilde{\sigma}^2 b(y)$}, with a scalar uncertainty parameter~\mbox{$\tilde{\sigma} > 0$}. In particular, it has been shown that, for $\tilde{\sigma} \rightarrow 0$, the tailored Jacobian approximation does not deteriorate the local convergence properties of the SQP algorithm. In the following, we translate these results to the case where~$b(y)$ is encoding uncertainty based on the GP posterior covariance by defining a parameter $\sigma := \sup_{y \in \mathcal{B}(\bar{y},r_\gamma)} \| b(y) \|$, which implicitly bounds the process noise and GP posterior covariance in a small neighborhood around a fixed point~$(\bar{y},\bar{P})$ of Alg.~\ref{alg:gpzoro}. This way, for GPs based on a twice continuously differentiable kernel, we can show that, for small enough process noise and GP posterior covariance, applying the tailored Jacobian approximation to GP-MPC preserves the local convergence properties of the SQP algorithm.

To this end, we first provide a brief review of the standard arguments for local
convergence analysis of inexact Newton-type methods using
strongly regular generalized equations.


\subsection{Generalized equations}

Generalized equations~\cite{robinsonStronglyRegularGeneralized1980} allow us to reformulate the KKT conditions of~\eqref{eq:OCP_compact} compactly as the set inclusion%
\begin{align}
    0 \in F(z) + \mathcal{N}_K(z), \label{eq:GE_optimal}
\end{align}
where%
\begin{align}
    F(z) := -\begin{bmatrix}
        \nabla_{(y,P)} \mathcal{L} (z) \\
        f(y) \\
        A(y) P + b(y) \\
        \bar{h}(y,P)
    \end{bmatrix} \label{eq:GE_optimal_F}
\end{align}
and $z := (y,P,\lambda^\mu,\lambda^\Sigma,\nu)$ contains the primal variables~\mbox{$y,P$} and the Lagrange multipliers~\mbox{$\lambda^\mu,\lambda^\Sigma,\nu$} corresponding to the constraints~\eqref{eq:OCP_compact_equality_mean}-\eqref{eq:OCP_compact_inequality}, respectively. Thereby,~$\mathcal{N}_K(z)$ denotes the normal cone to the set~\mbox{$K := \mathbb{R}^{n_y+n_P} \times \mathbb{R}^{n_f} \times \mathbb{R}^{n_g} \times \mathbb{R}^{n_h}_+$} at~$z$. 

While solutions to~\eqref{eq:GE_optimal} are KKT points of the original problem~\eqref{eq:OCP_GPMPC}, the suboptimal solution obtained by Alg.~\ref{alg:gpzoro} will instead satisfy the perturbed generalized equation%
\begin{align}
    0 \in \tilde{F}(z) + \mathcal{N}_K(z), \label{eq:GE_perturbed}
\end{align}
where%
\begin{align}
    \tilde{F}(z) := -\begin{bmatrix}
        \nabla_{(y,P)} \tilde{\mathcal{L}} (z) \\
        f(y) \\
        A(y) P + b(y) \\
        \bar{h}(y,P)
    \end{bmatrix}. \label{eq:GE_perturbed_F}
\end{align}
Evidently, the suboptimality is thereby caused by the inexact Jacobian approximation used in the QP subproblems, leading to a perturbed Lagrangian gradient%
\begin{align}
    \nabla_{(y,P)} \tilde{\mathcal{L}} (z) = \nabla_{(y,P)} \mathcal{L} (z) - 
    \begin{bmatrix}
        \frac{\partial}{\partial y} \left( A(y) P + b(y) \right)^\top \\
        0
    \end{bmatrix} \lambda^\Sigma
\end{align}
in the stationarity conditions in~\eqref{eq:GE_optimal_F} and~\eqref{eq:GE_perturbed_F}. 
Every iterate of Alg.~\ref{alg:gpzoro} solves the linearized generalized equation%
\begin{align}
    0 \in \tilde{F}(\hat{z}) + J(\hat{z})(z- \hat{z}) + \mathcal{N}_K(z), \label{eq:GE_zoro_lin}
\end{align}
where $J(\hat{z}) \approx \left. \frac{\partial \tilde{F}}{\partial z} \right|_{z=\hat{z}}$ is the Jacobian approximation around the linearization point $\hat{z}$. 

To obtain a sufficient local convergence criterion for Alg.~\ref{alg:gpzoro} to a solution~\mbox{$\bar{z} := (\bar{y}, \bar{P}, \bar{\lambda}^\mu, \bar{\lambda}^\Sigma, \bar{\nu})$} of~\eqref{eq:GE_perturbed}, we need the following key requirements to be satisfied. Thereby, $\mathcal{B}(z,r)$ denotes a ball of radius~\mbox{$r \in \mathbb{R}^{}$} centered at~\mbox{$z \in \mathbb{R}^{n_z}$}. 

\begin{assumption}[cf. 
    \cite{zanelliContractionEstimatesAbstract2019}, Ass.~1] 
    \label{ass:strong_regularity}
    Let~\eqref{eq:GE_perturbed} be strongly regular\footnote{Roughly speaking, the strong regularity assumption implies that there exists a single-valued and Lipschitz-continuous localization of the solution map (as defined in~\cite[p.4]{dontchevImplicitFunctionsSolution2009}) of~\eqref{eq:GE_perturbed} with respect to small perturbations around a point $\bar{z}$, which can also be stated in terms of a nonsingularity condition on the Jacobian of $\tilde{F}(z)$, c.f. the implicit-function theorem in the case of a fixed active set~\cite[Thm.~1B.1]{dontchevImplicitFunctionsSolution2009}. As such, by assuming strong regularity, we implicitly also assume differentiability of $\tilde{F}(z)$ at $\bar{z}$. For the NLP~\eqref{eq:OCP_compact} in particular, this requires differentiability of the square-root terms present in the tightened constraints~\eqref{eq:OCP_compact_inequality} at the solution, achievable by, e.g., a suitable design of $h(y_i)$, $B$ and $\Sigma^w$.} at $\bar{z}$, with Lipschitz constant~$\gamma$ in the neighborhood~$\mathcal{B}(\bar{z},r_\gamma)$ defined by $r_\gamma > 0$.
\end{assumption}

\begin{assumption}[cf. 
    \cite{zanelliContractionEstimatesAbstract2019}, Ass.~3]
    \label{ass:jacobian_diff}
    Let there exist a neighborhood $\mathcal{B}(\bar{z},r_{\tilde{\kappa}})$, with $0 < r_{\tilde{\kappa}} < r_\gamma$, and a positive constant $\tilde{\kappa}$, with $\gamma \tilde{\kappa} < \frac{1}{2}$, such that, for any $\hat{z} \in \mathcal{B}(\bar{z},r_{\tilde{\kappa}})$, it holds that%
    \begin{align}
        \left\| J(\hat{z}) - \left. \frac{\partial \tilde{F}}{\partial z} \right|_{z = \bar{z}} \right\| \leq \tilde{\kappa}.
    \end{align}
\end{assumption}
\vspace{1ex}
Under Assumptions~\ref{ass:strong_regularity} and~\ref{ass:jacobian_diff}, a sufficient local convergence criterion for Alg.~\ref{alg:gpzoro} reads as follows.

\begin{lemma}[cf.~\cite{zanelliContractionEstimatesAbstract2019},~Lemma~2]
    \label{lemma:GE_convergence}
    Let Assumptions~\ref{ass:strong_regularity} and \ref{ass:jacobian_diff} hold. Denote by $z_+$ a solution to~\eqref{eq:GE_zoro_lin} constructed at the linearization point~$\hat{z}$. Then, there exist strictly positive constants $\kappa < 1$ and  $r_\kappa$, such that, for any $\hat{z} \in \mathcal{B}(\bar{z},r_\kappa)$, it holds that%
    \begin{align}
        \| z_+ - \bar{z} \| \leq \kappa \| \hat{z} - \bar{z} \|. 
    \end{align} 
\end{lemma}

\vspace{1ex}

Note that for Newton-type optimization, Assumption~\ref{ass:strong_regularity} is standard and could be replaced by the stronger, but more frequently encountered, assumptions of linear independence constraint qualification and strong second-order sufficient condition, see~\cite[Thm.~4.1]{robinsonStronglyRegularGeneralized1980}. Regarding Assumption~\ref{ass:jacobian_diff}, the goal for the following section is to show that it is satisfied for sufficiently small process noise and GP posterior covariance matrices, ensuring local convergence of Alg.~\ref{alg:gpzoro} by Lemma~\ref{lemma:GE_convergence}.

\subsection{Local convergence for small uncertainties} \label{sec:local_convergence_proof}


The difference between the approximate and exact Jacobian of~\eqref{eq:GE_perturbed_F}, at the linearization point~$\hat{z}$ and a solution~$\bar{z}$ of~\eqref{eq:GE_perturbed}, respectively, is given by%
\begin{align}
    J(\hat{z}) - \left. \frac{\partial \tilde{F}}{\partial z} \right|_{z = \bar{z}} = \begin{bmatrix}
        M_{yy} - \left. \frac{\partial^2 \tilde{\mathcal{L}}}{\partial y^2} \right|_{z=\bar{z}} & 0 \\
        0 & 0 \\
        \left. \frac{\partial}{\partial y} \left( A(y) P + b(y) \right) \right|_{\overset{\scriptstyle y=\bar{y}}{P=\bar{P}}} & 0 \\
        0 & 0
    \end{bmatrix}. \label{eq:GE_jacobian_diff}
\end{align}
Hence, from Lemma~\ref{lemma:GE_convergence} and equation~\eqref{eq:GE_jacobian_diff}, we can deduce local convergence in a neighborhood around~$\bar{z}$ if the errors in the Hessian- and tailored Jacobian approximations are sufficiently small. 
As only the latter is essential to our method, we focus on the Jacobian difference induced by the tailored Jacobian approximation; see e.g.~\cite{boggsSequentialQuadraticProgramming1995} for a discussion of the Hessian approximation's role in SQP methods. 

The main result makes use of the following regularity assumption. 
\begin{assumption} \label{ass:b_diffable}
    Let $b(y)$ be twice continuously differentiable for all $y \in \mathcal{B}(\bar{y}, r_\gamma)$.
\end{assumption}
Note that this assumption can also be phrased in terms of the GP kernel function's regularity: As~$b(y)$ denotes the vectorized process noise and GP posterior covariances, the latter of which are a linear combination of kernel function evaluations~\cite{rasmussenGaussianProcessesMachine2006}, we restrict ourselves to twice continuously differentiable kernels. This captures many of the kernels commonly used in practice, such as squared exponential, linear or Matèrn kernels with $\nu \geq 5/2$. 

\begin{lemma} \label{lemma:jacobian_diff_bound}
    Let Assumptions~\ref{ass:strong_regularity} and~\ref{ass:b_diffable} hold and define \mbox{$\sigma := \sup_{y \in \mathcal{B}(\bar{y}, r_\gamma)} {\| b(y) \|}$}.
    Then, for any $\epsilon > 0$, there exists a $\delta \in \mathbb{R}^{}$ such that if $\sigma \leq \delta$, it holds that%
    \begin{align}
        \left\|  \left. \frac{\partial}{\partial y} \left( A(y) P + b(y) \right) \right|_{\overset{\scriptstyle y = \bar{y}}{P = \bar{P}}} \right\| \leq \epsilon. \label{eq:GE_jacobian_constraint_diff}
    \end{align}

    \begin{proof}
        We will show the above implication by considering both summands inside~\eqref{eq:GE_jacobian_constraint_diff} separately and applying the triangle inequality. 
        
        For the first term, rearranging and taking the norm of the equality constraint~\eqref{eq:prop_linear}, evaluated at the solution $\bar{z}$, leads to%
        \begin{align*}
            \| \bar{P} \| = \| A(\bar{y})^{-1} b(\bar{y}) \| \leq \| A(\bar{y})^{-1} \| \| b(\bar{y}) \| \leq \| A(\bar{y})^{-1} \| \sigma.
        \end{align*}
        Hence, 
        by linearity of%
        ~\mbox{$\frac{\partial}{\partial y} A(y) P$}
        in $P$ and the triangle inequality, 
        for any \mbox{$\epsilon > 0$}, there exists a~\mbox{$\delta \in \mathbb{R}^{}$} such that $\sigma \leq \delta$ implies \mbox{$\left\| \left. \frac{\partial}{\partial y} A(y) \bar{P} \right|_{y = \bar{y}} \right\| \leq \epsilon$}. 

        For the second term, 
        consider the Taylor expansion of the $i$-th component of $b(y)$ around the expansion point $\bar{y}$, evaluated at $y_t \in \mathcal{B}(\bar{y},r_\gamma)$,
        \begin{align}
            b_i(y_t) = b_i(\bar{y}) + t \left. \frac{\partial b_i(y)}{\partial y} \right|_{y=\bar{y}} \delta y + \frac{t^2}{2} \delta y^\top \left. \frac{\partial^2 b_i(y)}{\partial y^2} \right|_{y=\xi_i} \delta y, \label{eq:GE_proof_taylor}
        \end{align}
        where $t := \| y_t - \bar{y} \|$ and $\delta y := (y_t - \bar{y})/t$ (for $y_t \neq \bar{y}$) is a unit vector. By the mean-value theorem, the Hessian is evaluated at some $\xi_i \in \Xi$ for each component~$i$, with the set
        \mbox{$\Xi := \left\{ \bar{y} + s (y_t - \bar{y}) \> \middle| \> s \in [0,1] \text{ and } \bar{y} \in \mathbb{R}^{n_y}, y_t \in \mathcal{B}(\bar{y},r_\gamma) \right\}$}. 
        Since $b_i(y)$ is twice continuously differentiable, there exists some $M_2 \in \mathbb{R}^{}$ such that, for all $\xi_i \in \Xi$, $i \in \{ 1,\ldots,n_x^2(N+1) \}$, it holds that%
        \begin{align}
            \left\| \left. \frac{\partial^2 b_i(y)}{\partial y^2} \right|_{y=\xi_i} \right\| \leq M_2.
        \end{align}
        From~\eqref{eq:GE_proof_taylor}, by taking the absolute value, dividing by $t > 0$ and inserting the definitions of $\sigma$ and $M_2$, we obtain%
        \begin{align}
            \left| \left. \frac{\partial b_i(y)}{\partial y} \right|_{y=\bar{y}} \delta y \right| &= \left| \frac{b_i(y_t) - b_i(\bar{y})}{t} - \frac{t}{2} \delta y^\top \left. \frac{\partial^2 b_i(y)}{\partial y^2} \right|_{y=\xi_i} \delta y \right| \notag \\
            &\leq \frac{2\sigma}{t} + \frac{M_2}{2} t. \label{eq:GE_proof_inequality_before_sup}
        \end{align}
        Since~\eqref{eq:GE_proof_taylor} holds for all $0 < t < r_\gamma$ and $\delta y$ with $\| \delta y \| = 1$, taking the supremum over $\delta y$ on the left and the infimum over $t$ on the right-hand side preserves inequality~\eqref{eq:GE_proof_inequality_before_sup}, i.e.,%
        \begin{align}
            \sup_{\| \delta y \| = 1} \left| \left. \frac{\partial b_i(y)}{\partial y} \right|_{y=\bar{y}} \delta y \right| &\leq \inf_{t>0} \frac{2\sigma}{t} + \frac{M_2}{2} t. \label{eq:GE_proof_inequality_infsup}
            \end{align}
            By applying the definition of the induced matrix norm on the left, and solving for the value of the infimum on the right-hand side ($t < r_\gamma$ can always be achieved by increasing $M_2$), \eqref{eq:GE_proof_inequality_infsup}~simplifies to
            \begin{align}
                \left\| \left. \frac{\partial b_i(y)}{\partial y} \right|_{y=\bar{y}} \right\| &\leq 2\sqrt{\sigma M_2}.
        \end{align}

        Thus, we have shown that the norms of both summands in~\eqref{eq:GE_jacobian_constraint_diff} scale with~\mbox{$\mathcal{O}(\sigma)$}, which proves the assertion by the triangle inequality and the definition of~$\mathcal{O}(\sigma)$.
    \end{proof}
\end{lemma}

To summarize, in tandem with Lemma~\ref{lemma:GE_convergence} and a sufficiently accurate Hessian approximation~$M_{yy}$ across all SQP iterates, Lemma~\ref{lemma:jacobian_diff_bound} establishes guaranteed convergence of Alg.~\ref{alg:gpzoro} to a solution~$\bar{z}$ of~\eqref{eq:GE_perturbed} for sufficiently small process noise and GP covariances in a local neighborhood and, practical convergence depending on the regularity of~\eqref{eq:GE_perturbed} in terms of the Lipschitz constant~$\gamma$. 




\section{Numerical results} \label{sec:experiments}

Alg.~\ref{alg:gpzoro} has been prototyped in \texttt{Python} using the corresponding \texttt{acados} interface. For the nominal dynamics, evaluation and sensitivity computation is thereby performed using an \texttt{acados} integrator and just-in-time compiled \texttt{CasADi} functions~\cite{anderssonCasADiSoftwareFramework2019}; for the GP~mean and covariance, the corresponding computations are carried out with \texttt{PyTorch}~\cite{paszkePyTorchImperativeStyle2019} using the \texttt{GPyTorch} library~\cite{gardnerGPyTorchBlackboxMatrixMatrix2018}.\footnote{An open-source implementation of the following example is available at \url{https://gitlab.ethz.ch/ics/zero-order-gp-mpc}, \doi{10.3929/ethz-b-000611298}.} 
In the following, we compare different variants of the proposed algorithm against a ``na\"ive'' GP-MPC implementation and nominal MPC using a scalable benchmarking example. 

\subsection{Hanging chain example}

To test the scaling properties of Alg.~\ref{alg:gpzoro}, we apply it on a slightly modified variant of the hanging chain example, a popular benchmark for numerical methods for NMPC~\cite{zanelliZeroOrderRobustNonlinear2021},~\cite{kouzoupisRecentAdvancesQuadratic2018}. 
The system is defined by a chain of 
masses \mbox{$m := 0.033\si{kg}$} connected by linear springs with stiffness \mbox{$k := 30.3\si{N/m}$}. The mass at one end of the chain is tied to the origin; the \mbox{$n_u = 3$} velocity components of the other end's mass constitute the control inputs of the system. The system state is given by the position components of the controlled mass as well as the position and velocity components of the intermediate masses, resulting in a state space dimension of \mbox{$n_x = 6(n_\text{mass} - 2) + 3$}. The initial state of the chain is computed based on its resting position, where the controlled end is placed at \mbox{$(x_{\text{init}},y_{\text{init}},z_{\text{init}}) := (6l (n_\text{mass} - 1), 0, 0)$}, with length~\mbox{$l := 0.033\si{m}$}. 
After applying a control input of \mbox{$u_{\text{init}} := (1, 1, 1)$} for \mbox{$T_{\text{init}} := 1\si{s}$}, the control task is to restore the resting position while ensuring that none of the masses violates a wall constraint at \mbox{$y_{\text{wall}}:= -0.05\si{m}$}. 
The example is modified by adding a latent force 
\begin{align*}
    f_\text{lat}(x,v_x) := \alpha_{\text{lat}} \left(v_x - \sin \left( \beta_1 \frac{2 \pi x}{l} \right) - \sin \left( \beta_2 \frac{2 \pi x}{l} \right)^2 \right)^2
\end{align*}
to the $y$-acceleration of each intermediate mass, where \mbox{$\alpha_{\text{lat}} := -0.1$}, \mbox{$\beta_1 := 2$}, \mbox{$\beta_1 := 3$}, and $x$ and $v_x$ denote the position and velocity along the \mbox{$x$-axis} of the frame, respectively. 
The continuous-time dynamics are discretized using an implicit Runge-Kutta integrator, with time step \mbox{$T_s := 0.2\si{s}$}. A Gauss-Newton Hessian approximation~$M_{yy}$ is employed. 
The model mismatch on each of the~\mbox{$n_w := 3(n_{\text{mass}} - 2)$} velocity states of the intermediate masses is captured using independent GPs with squared exponential kernel. 
Training data is generated by recording~\mbox{$D := 15 N_{x_0}$} samples of the model mismatch from closed-loop simulations with a nominal model predictive controller that only considers the nominal dynamics and no constraint tightening, starting from $N_{x_0}$ perturbed initial conditions. The following experiments were performed using an Intel~i9-7940X~processor running at 3.10~GHz and an NVIDIA~GeForce~RTX~2080~Ti~GPU. 

\subsection{Timings for increasing state and GP output dimension}

\begin{figure}    
    \centering
    \includegraphics[scale=1]{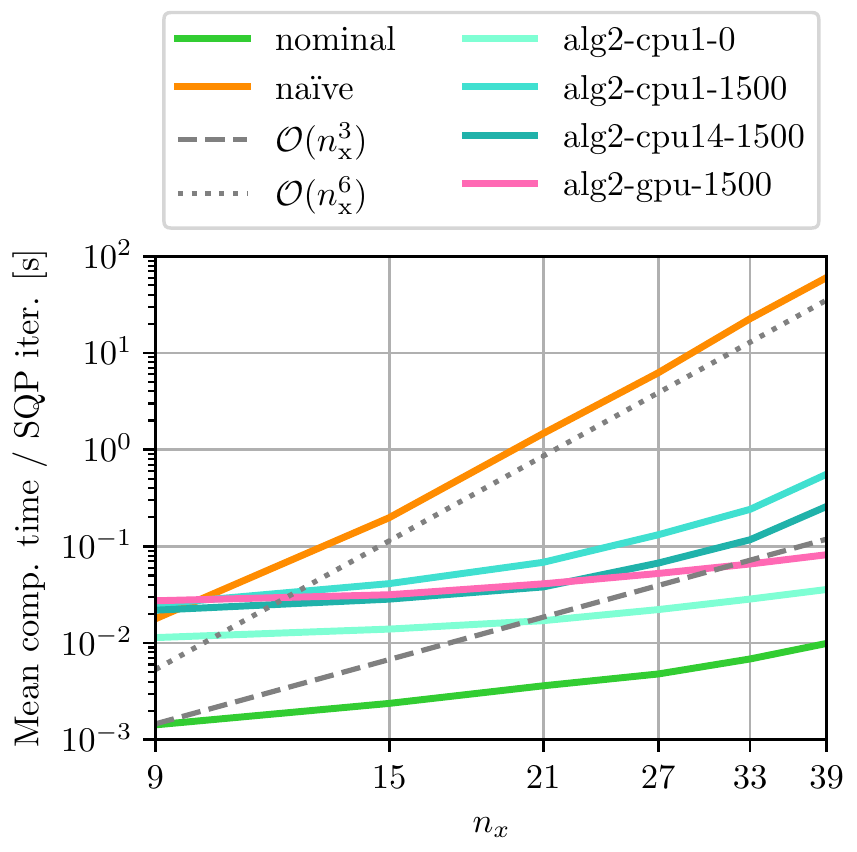}
    \caption{SQP timings comparison for increasing number of states $n_x$
    and GP dimension $n_w = (n_x-3)/2$. 
    }
    \label{fig:ChainTimings}
\end{figure}

Fig.~\ref{fig:ChainTimings} shows the scaling of the mean computation times per SQP iteration as the number of states~$n_x$ is increased by adding masses to the chain, for the following controller implementations:
\begin{itemize}
    \item ``nominal'': Using only the nominal model and no uncertainty description,
    \item ``na\"ive'': Solving~\eqref{eq:OCP_GPMPC} exactly without any data points (no GP-related computations\footnote{
        Due to software limitations, for a fair comparison the ``na\"ive'' \mbox{GP-MPC} implementation only makes use of the GP prior, i.e., \mbox{$\mu^d(y) \equiv 0$} and \mbox{$\Sigma^d(y) \equiv const.$}; it corresponds to a stochastic NMPC implementation with the process noise covariance inflated by the GP prior covariance. Therefore the present timings can be viewed as a lower bound for the actual timings when solving~\eqref{eq:OCP_GPMPC} with a GP conditioned on data, which would additionally require not only GP~posterior evaluations and derivative computations, but also expensive computations associated with the Hessian of the GP~posterior mean.
        }), including the covariances as optimization variables into an augmented state, 
    \item ``alg2-cpu$C$-$D$'': Alg.~\ref{alg:gpzoro} with $D$ data points and GP inference and AD on $C$ CPU cores, and
    \item ``alg2-gpu-$D$'': Alg.~\ref{alg:gpzoro} with $D$ data points and GP inference and AD performed on GPU. 
\end{itemize}

It becomes evident that 
the scaling of order $\mathcal{O}(n_x^6)$ for the ``na\"ive'' method quickly becomes prohibitive in terms of computation time. In comparison, for $n_x = 39$, a roughly $1000$-fold speed-up by means of the zero-order optimization strategy can be observed. 
Comparing the computation times of Alg.~\ref{alg:gpzoro} for different number of data points clearly shows the effects of GPU acceleration, which accounts for a $5-10$-fold speed-up and better scaling properties per SQP iteration for~\mbox{$n_x \geq 33$} and the associated GP~output dimension. 
                

\subsection{Timing profile of Alg.~\ref{alg:gpzoro}}

\begin{figure}
    \vspace{2mm}
    \centering
    \includegraphics[scale=1]{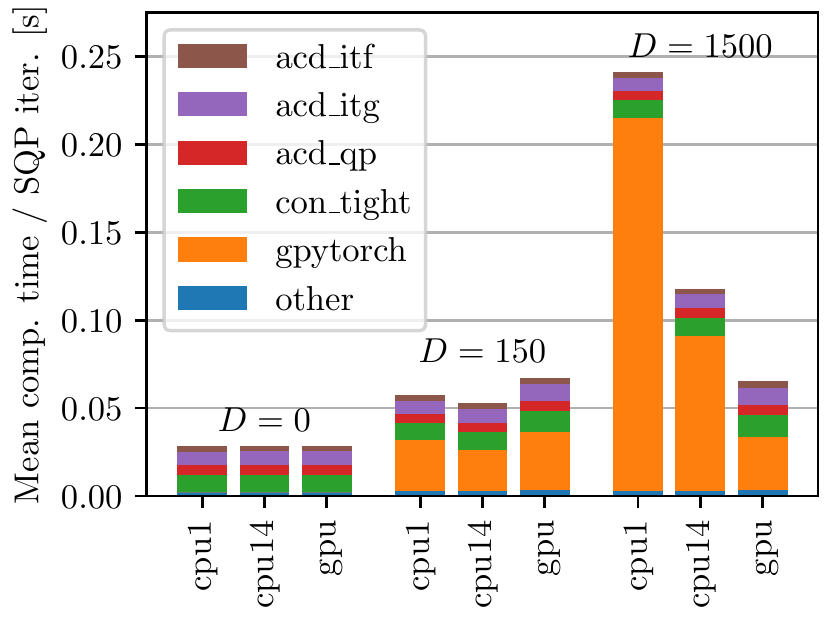}
    \caption{Timing profile for Alg.~\ref{alg:gpzoro} variants for $n_x = 33$ ($n_\text{mass} = 7$). 
    }
    \label{fig:RunTimesGPU}    
\end{figure}

Fig.~\ref{fig:RunTimesGPU} displays the runtime profile of Alg.~\ref{alg:gpzoro} for $n_x = 33$ states, divided into the following components:
\begin{itemize}
    \item ``acd\_itf'': \texttt{acados} interface for get/set operations,
    \item ``acd\_itg'': \texttt{acados} implicit Runge-Kutta integrator,
    \item ``acd\_qp'': \texttt{acados} QP solve of~\eqref{eq:QP_zoro} using HPIPM~\cite{frisonHPIPMHighperformanceQuadratic2020},
    \item ``con\_tight'': \texttt{NumPy} covariance propagation~\eqref{eq:CovariancePropagation} and constraint tightening for~\eqref{eq:QP_zoro_ineq},
    \item ``gpytorch'': \texttt{GPyTorch} GP inference and AD.
\end{itemize}
For 150 data points, Fig.~\ref{fig:RunTimesGPU} shows that the overhead introduced by parallelizing \texttt{gpytorch} operations is higher than the computational speed-ups, especially on the GPU. For 1500 data points, however, the parallelization leads to significantly lower computation times, saving about 60-70\% of the computational costs associated with the GPs, which at this point dominate the computational footprint of the method. For 0~data points, potential for improvement can be seen in terms of the constraint tightening computations, 
which could be improved by switching 
to a \texttt{C/C++} implementation.

\addtolength{\textheight}{-3.5cm}

\section{Conclusions}

To make GP-MPC computationally tractable, both the complexity to solve the OCP, as well as evaluating and differentiating the GP~posterior mean and covariance need to be addressed. To tackle these challenges, this paper has presented an inexact SQP approach with a tailored Jacobian approximation, while parallelizing GP inference and differentiation routines on a GPU. The results demonstrate that drastically reduced computation times can be achieved while ensuring feasibility of the converged iterates and maintaining favorable convergence properties, pushing the computational boundaries to apply GP-MPC in real-world scenarios. 

\subsection*{Acknowledgments}

We would like to thank Katrin Baumgärtner and Jonathan Frey for pointing us towards the \texttt{acados} \texttt{Cython} interface, and Johannes Köhler for many helpful discussions.


\bibliography{IEEEabrv,Zero-Order-GP-MPC}

\begin{thebibliography}{10}
\providecommand{\url}[1]{#1}
\csname url@rmstyle\endcsname
\providecommand{\newblock}{\relax}
\providecommand{\bibinfo}[2]{#2}
\providecommand\BIBentrySTDinterwordspacing{\spaceskip=0pt\relax}
\providecommand\BIBentryALTinterwordstretchfactor{4}
\providecommand\BIBentryALTinterwordspacing{\spaceskip=\fontdimen2\font plus
\BIBentryALTinterwordstretchfactor\fontdimen3\font minus
  \fontdimen4\font\relax}
\providecommand\BIBforeignlanguage[2]{{%
\expandafter\ifx\csname l@#1\endcsname\relax
\typeout{** WARNING: IEEEtran.bst: No hyphenation pattern has been}%
\typeout{** loaded for the language `#1'. Using the pattern for}%
\typeout{** the default language instead.}%
\else
\language=\csname l@#1\endcsname
\fi
#2}}

\bibitem{ostafewLearningbasedNonlinearModel2014}
C.~J. Ostafew, A.~P. Schoellig, and T.~D. Barfoot, ``Learning-based nonlinear
  model predictive control to improve vision-based mobile robot path-tracking
  in challenging outdoor environments,'' in \emph{2014 {{IEEE International
  Conference}} on {{Robotics}} and {{Automation}} ({{ICRA}})}, 2014, pp.
  4029--4036.

\bibitem{carronDataDrivenModelPredictive2019}
A.~Carron, E.~Arcari, M.~Wermelinger, L.~Hewing, M.~Hutter, and M.~N.
  Zeilinger, ``Data-{{Driven Model Predictive Control}} for {{Trajectory
  Tracking With}} a {{Robotic Arm}},'' \emph{IEEE Robotics and Automation
  Letters}, vol.~4, no.~4, pp. 3758--3765, 2019.

\bibitem{kabzanLearningBasedModelPredictive2019}
J.~Kabzan, L.~Hewing, A.~Liniger, and M.~N. Zeilinger, ``Learning-{{Based Model
  Predictive Control}} for {{Autonomous Racing}},'' \emph{IEEE Robot. Autom.
  Lett.}, vol.~4, no.~4, pp. 3363--3370, 2019.

\bibitem{kabzanAMZDriverlessFull2020}
J.~Kabzan, M.~Valls, V.~Reijgwart, H.~Hendrikx, C.~Ehmke, M.~Prajapat,
  A.~B{\"u}hler, N.~Gosala, M.~Gupta, R.~Sivanesan, A.~Dhall, E.~Chisari,
  N.~Karnchanachari, S.~Brits, M.~Dangel, I.~Sa, R.~Dube, A.~Gawel,
  M.~Pfeiffer, and R.~Siegwart, ``{{AMZ Driverless}}: {{The}} full autonomous
  racing system,'' \emph{Journal of Field Robotics}, 2020.

\bibitem{torrenteDataDrivenMPCQuadrotors2021}
G.~Torrente, E.~Kaufmann, P.~F{\"o}hn, and D.~Scaramuzza, ``Data-{{Driven MPC}}
  for {{Quadrotors}},'' \emph{IEEE Robotics and Automation Letters}, vol.~6,
  no.~2, pp. 3769--3776, 2021.

\bibitem{hewingCautiousModelPredictive2020}
L.~Hewing, J.~Kabzan, and M.~N. Zeilinger, ``Cautious {{Model Predictive
  Control Using Gaussian Process Regression}},'' \emph{IEEE Trans Control Syst
  Technol}, vol.~28, no.~6, pp. 2736--2743, 2020.

\bibitem{frohlichModelLearningContextual2022}
L.~P. Fr{\"o}hlich, C.~K{\"u}ttel, E.~Arcari, L.~Hewing, M.~N. Zeilinger, and
  A.~Carron, ``Model {{Learning}} and {{Contextual Controller Tuning}} for
  {{Autonomous Racing}},'' \emph{2022 IEEE/RSJ International Conference on
  Intelligent Robots and Systems (IROS), to appear}, 2022.

\bibitem{vaskovFrictionAdaptiveStochasticPredictive2022}
S.~Vaskov, R.~Quirynen, M.~Menner, and K.~Berntorp, ``Friction-{{Adaptive
  Stochastic Predictive Control}} for {{Trajectory Tracking}} of {{Autonomous
  Vehicles}},'' in \emph{2022 {{American Control Conference}} ({{ACC}})}, 2022,
  pp. 1970--1975.

\bibitem{fengInexactAdjointbasedSQP2020}
X.~Feng, S.~D. Cairano, and R.~Quirynen, ``Inexact {{Adjoint-based SQP
  Algorithm}} for {{Real-Time Stochastic Nonlinear MPC}},''
  \emph{IFAC-PapersOnLine}, vol.~53, no.~2, pp. 6529--6535, 2020.

\bibitem{zanelliZeroOrderRobustNonlinear2021}
A.~Zanelli, J.~Frey, F.~Messerer, and M.~Diehl, ``Zero-{{Order Robust Nonlinear
  Model Predictive Control}} with {{Ellipsoidal Uncertainty Sets}},''
  \emph{IFAC-PapersOnLine}, vol.~54, no.~6, pp. 50--57, 2021.

\bibitem{paszkePyTorchImperativeStyle2019}
A.~Paszke, S.~Gross, F.~Massa, A.~Lerer, J.~Bradbury, G.~Chanan, T.~Killeen,
  Z.~Lin, N.~Gimelshein, L.~Antiga, A.~Desmaison, A.~Kopf, E.~Yang, Z.~DeVito,
  M.~Raison, A.~Tejani, S.~Chilamkurthy, B.~Steiner, L.~Fang, J.~Bai, and
  S.~Chintala, ``{{PyTorch}}: {{An Imperative Style}}, {{High-Performance Deep
  Learning Library}},'' in \emph{Advances in {{Neural Information Processing
  Systems}} 32}.\hskip 1em plus 0.5em minus 0.4em\relax {Curran Associates,
  Inc.}, 2019, pp. 8024--8035.

\bibitem{gardnerGPyTorchBlackboxMatrixMatrix2018}
J.~Gardner, G.~Pleiss, K.~Q. Weinberger, D.~Bindel, and A.~G. Wilson,
  ``{{GPyTorch}}: {{Blackbox Matrix-Matrix Gaussian Process Inference}} with
  {{GPU Acceleration}},'' in \emph{Advances in {{Neural Information Processing
  Systems}}}, vol.~31.\hskip 1em plus 0.5em minus 0.4em\relax {Curran
  Associates, Inc.}, 2018.

\bibitem{messererEfficientAlgorithmTubebased2021}
F.~Messerer and M.~Diehl, ``An {{Efficient Algorithm}} for {{Tube-based Robust
  Nonlinear Optimal Control}} with {{Optimal Linear Feedback}},'' in \emph{2021
  60th {{IEEE Conference}} on {{Decision}} and {{Control}} ({{CDC}})}, 2021,
  pp. 6714--6721.

\bibitem{quirynenUncertaintyPropagationLinear2021}
R.~Quirynen and K.~Berntorp, ``Uncertainty {{Propagation}} by {{Linear
  Regression Kalman Filters}} for {{Stochastic NMPC}},''
  \emph{IFAC-PapersOnLine}, vol.~54, no.~6, pp. 76--82, 2021.

\bibitem{rasmussenGaussianProcessesMachine2006}
C.~E. Rasmussen and C.~K.~I. Williams, \emph{Gaussian Processes for Machine
  Learning}, ser. Adaptive Computation and Machine Learning.\hskip 1em plus
  0.5em minus 0.4em\relax {Cambridge, Massachusetts}: {MIT Press}, 2006.

\bibitem{girardGaussianProcessPriors2002}
A.~Girard, C.~E. Rasmussen, and R.~{Murray-Smith}, ``Gaussian {{Process}}
  priors with {{Uncertain Inputs}}: {{Multiple-Step-Ahead Prediction}},''
  {Department of Computing Science, University of Glasgow}, Technical
  {{Report}} TR-2002-119, 2002.

\bibitem{hewingCautiousNMPCGaussian2018}
L.~Hewing, A.~Liniger, and M.~N. Zeilinger, ``Cautious {{NMPC}} with {{Gaussian
  Process Dynamics}} for {{Autonomous Miniature Race Cars}},'' in \emph{2018
  {{European Control Conference}} ({{ECC}})}, 2018, pp. 1341--1348.

\bibitem{bockMultipleShootingAlgorithm1984}
H.~G. Bock and K.~J. Plitt, ``A {{Multiple Shooting Algorithm}} for {{Direct
  Solution}} of {{Optimal Control Problems}}*,'' \emph{IFAC Proceedings
  Volumes}, vol.~17, no.~2, pp. 1603--1608, 1984.

\bibitem{kouzoupisRecentAdvancesQuadratic2018}
D.~Kouzoupis, G.~Frison, A.~Zanelli, and M.~Diehl, ``Recent {{Advances}} in
  {{Quadratic Programming Algorithms}} for {{Nonlinear Model Predictive
  Control}},'' \emph{Vietnam J. Math.}, vol.~46, no.~4, pp. 863--882, 2018.

\bibitem{boggsSequentialQuadraticProgramming1995}
P.~T. Boggs and J.~W. Tolle, ``Sequential quadratic programming,'' \emph{Acta
  Numerica}, vol.~4, pp. 1--51, 1995.

\bibitem{bockConstrainedOptimalFeedback2007}
H.~G. Bock, M.~Diehl, E.~Kostina, and J.~P. Schl{\"o}der, ``1. {{Constrained
  Optimal Feedback Control}} of {{Systems Governed}} by {{Large Differential
  Algebraic Equations}},'' in \emph{Real-{{Time PDE-Constrained
  Optimization}}}.\hskip 1em plus 0.5em minus 0.4em\relax {Society for
  Industrial and Applied Mathematics}, 2007, pp. 3--24.

\bibitem{robinsonStronglyRegularGeneralized1980}
S.~M. Robinson, ``Strongly {{Regular Generalized Equations}},'' \emph{Math.
  Oper. Res.}, vol.~5, no.~1, pp. 43--62, 1980.

\bibitem{zanelliContractionEstimatesAbstract2019}
A.~Zanelli, Q.~{Tran-Dinh}, and M.~Diehl, ``Contraction {{Estimates}} for
  {{Abstract Real-Time Algorithms}} for {{NMPC}},'' in \emph{2019 {{IEEE}} 58th
  {{Conference}} on {{Decision}} and {{Control}} ({{CDC}})}, 2019, pp.
  8085--8092.

\bibitem{dontchevImplicitFunctionsSolution2009}
A.~L. Dontchev and R.~T. Rockafellar, \emph{Implicit {{Functions}} and
  {{Solution Mappings}}: {{A View}} from {{Variational Analysis}}}, ser.
  Springer {{Monographs}} in {{Mathematics}}.\hskip 1em plus 0.5em minus
  0.4em\relax {New York, NY}: {Springer}, 2009.

\bibitem{anderssonCasADiSoftwareFramework2019}
J.~A.~E. Andersson, J.~Gillis, G.~Horn, J.~B. Rawlings, and M.~Diehl,
  ``{{CasADi}}: A software framework for nonlinear optimization and optimal
  control,'' \emph{Math. Prog. Comp.}, vol.~11, no.~1, pp. 1--36, 2019.

\bibitem{frisonHPIPMHighperformanceQuadratic2020}
G.~Frison and M.~Diehl, ``{{HPIPM}}: A high-performance quadratic programming
  framework for model predictive control,'' \emph{IFAC-PapersOnLine}, vol.~53,
  no.~2, pp. 6563--6569, 2020.

\end{thebibliography}

\appendix

\subsection{Explicit form of~\eqref{eq:prop_linear}} \label{sec:appendix_matrices}

By linearity of $g(y,P)$ in $P$, we can show that $\frac{\partial g}{\partial P}(y,P)$ is indeed invertible for our problem setting. 
It holds that
\begin{align}
    A(y) &:= \begin{bmatrix}
        I_{n_x^2} & &\\
        L(y_0) & I_{n_x^2} & \\
        & \ddots & \ddots & \\ 
        & & L(y_{N-1}) & I_{n_x^2} \\ 
    \end{bmatrix}
    \intertext{and}
    b(y) &:= \begin{bmatrix}
        0 \\
        (B \otimes B) \overrightarrow{\left( \Sigma^d(y_0) + \Sigma^w \right)} \\
        \vdots \\
        (B \otimes B) \overrightarrow{\left( \Sigma^d(y_{N-1}) + \Sigma^w \right)}
    \end{bmatrix},
\end{align}
where $L(y_i) := -\tilde{A}_i(y_i) \otimes \tilde{A}_i(y_i)$ and ``$\otimes$'' denotes the Kronecker product. Differentiating equation~\eqref{eq:prop_linear}, we obtain that \mbox{$\frac{\partial g}{\partial P}(y) = A(y)$}; invertibility of $\frac{\partial g}{\partial P}$ follows since $A(y)$ is lower triangular. 

\end{document}